\documentclass{amsart} 
\usepackage{amsfonts,amsmath,amsthm,amssymb}
\usepackage{clrscode}
\usepackage{url}
\begin{document}

\title{{\Small CIRCA preprint 2010/9\vspace{10pt}\\} Rewriting the check of 8-rewritability for $A_5$}

\author{Alexander Konovalov}
\address{School of Computer Science, University of St Andrews, St Andrews, Fife, KY16 9SX Scotland}
\email{alexk@mcs.st-andrews.ac.uk}




\begin{abstract}
The group $G$ is called $n$-rewritable for $n>1$, 
if for each sequence of $n$ elements $x_1, x_2, \dots, x_n \in G$
there exists a non-identity permutation $\sigma \in S_n$ such that
$x_1 x_2 \cdots x_n = x_{\sigma(1)} x_{\sigma(2)} \cdots x_{\sigma(n)}$.
Using computers, Blyth and Robinson (1990) verified that 
the alternating group $A_5$ is 8-rewritable. We report on an 
independent verification of this statement 
using the computational algebra system GAP, and compare the
performance of our sequential and parallel code with the original 
one.
\end{abstract}

\maketitle

Let $n>1$ be an integer. Following \cite{Blyth-Robinson},
a group $G$ is said to be \emph{totally $n$-rewritable},
or have the rewriting property $P_n$, if for each sequence
of $n$ elements $x_1, x_2, \dots, x_n$ of the group $G$ there
exists a non-identity permutation $\sigma \in S_n$ such that
$$
x_1 x_2 \cdots x_n = x_{\sigma(1)} x_{\sigma(2)} \cdots x_{\sigma(n)}.
$$
Clearly, all abelian groups satisfy $P_2$, and if $G$ satisfies
$P_{k}$ then it also satisfies $P_{k+1}$.

On the problem session of the conference
``Arithmetic of Group Rings and Related Objects'' 
(Aachen, Germany, March 22-26, 2010) 
Eli Aljadeff (Technion, Haifa, Israel) 
suggested the following problem:\\

\centerline{
Prove that the alternating group $A_5$ 
has the property $P_8$.} 
\vspace{10pt}

\noindent
He referred to the computer verification of this statement 
reported in \cite{Blyth-Robinson}, and demonstrated how to show
that $A_5$ has the property $P_{10}$ using group rings technique. 
He also suggested that since \cite{Blyth-Robinson} appeared 
twenty years ago, nowadays this result probably could be verified
much faster. Motivated by this, the author verified that $A_5$ has 
the property $P_8$ using the computational algebra system GAP \cite{GAP}
and compared the performance of the sequential and parallel GAP 
implementations with the one described in \cite{Blyth-Robinson}.

To check that the group $G$ is $n$-rewritable using the brute force
approach, one may enumerate all $n$-tuples of distinct elements
of $G$ and check that each of them may be rewritten. Of course, even 
for $A_5$ the number of tuples to check will be enormous, so this
approach will not work.

There is, however, a simple observation that allows to reduce the number
of checks substantially. The algorithm described by Blyth and Robinson 
in \cite{Blyth-Robinson} constructs all non-rewritable words of length 2
which are non-equivalent with respect to the action of ${\rm Aut}(G)$. 
On the next step, these words are
used to construct all non-equivalent non-rewritable words of length 3,
and then the process is repeated until the first $n$ for which there are
no non-rewritable words of length $n$ will be found. 

The next table contains information about the number $N(r)$ of 
non-rewritable words of length $r$ in $A_5$, determined using computers 
and listed in \cite{Blyth-Robinson}:
$$
\begin{array}{|c|l|}
\hline
r & N(r) \\
\hline
2 & 29 \\
3 & 1315 \\
4 & 43121 \\
5 & 528069 \\
6 & 187719 \\
7 & 1320 \\
8 & 0 \\
\hline
\end{array}
$$
The authors of \cite{Blyth-Robinson} wrote that first these data were computed over
a period of two weeks by a PASCAL program on a MicroVAX II, and then verified by 
a parallel C++ implementation that produced the same result on four Sun 3/60 machines
in less than three hours. 

We were interested to compare the reported performance with the runtime that can be
achieved nowadays in a sequential version of GAP on modern computers. To test our
implementation, we used an 8-core Intel server, with dual quad-core Intel Xeon 5570 
2.93GHz / RAM 48 GB / CentOS Linux 5.3.
First we present the algorithm in the pseudocode, following its
textual description from \cite{Blyth-Robinson}:

\vspace{20pt}

\begin{codebox}
\Procname{$\proc{RewritabilityLength}(G,m)$}
\li $A \gets \proc{AutomorphismGroup}[G]$
\li $x \gets \proc{NontrivialOrbitRepresentatives}[A,G]$
\li $n \gets 1$
\li \Repeat
\li 	$n \gets n+1$
\li 	$\id{nrw} \gets \const{empty list}$
\li 	\For $u$ \kw{in} $x$ \>\>\>\>\Comment $u$ is a non-rewritable word of length $n-1$
\li 	\Do $K \gets \proc{Intersection}[$
	                 \Indentmore
\zi 	             $\proc{Stabiliser}[A,u[1]],$
\zi 	             $ \dots,$ 
\zi 	             $\proc{Stabiliser}[A,u[n-1]]]$
	                 \End
\li 		\If $\proc{Size}[K] > 1$
\li 			\Then $y \gets \proc{NontrivialOrbitRepresentatives}[K,G]$
\li 			\Else $y \gets G \setminus \{ 1_G \}$
			\End
\li 		\For $v$ \kw{in} $y$
\li 		\Do $t = \proc{Concatenation}[u,v]$
\li             \If $\proc{IsRewritableWord}[t]$
\li 	        \Then $\proc{Append}[nrw,t]$
	            \End
	        \End   
	    \End     
\li 	\If $\id{nrw} = \const{empty list}$
\li 	\Then \Return $n$
\li 	\Else $x \gets \id{nrw}$
	\End
\li \Until $n = m$
\li \Return \id{fail}
\End
\end{codebox}

\vspace{20pt}

The pseudocode above refers to the following procedures: 
\begin{itemize}
\item $\proc{AutomorphismGroup}[G]$ returns the automorphism group of $G$;
\item $\proc{NontrivialOrbitRepresentatives}[A,G]$ return the list of representatives
      of orbits of non-identity elements of the group $G$ under the action of its automorphism
      group $A$;
\item $\proc{Stabiliser}[A,g]$ returns the stabiliser of an element $g \in G$ in the group $A$;
\item $\proc{IsRewritableWord}[t]$ checks if the word $t$ is rewritable.
\end{itemize}
Other names of procedures should be self-explanatory. 

The first implementation looked very much like the pseudocode above, 
and it took almost 34 hours to run (though it used one CPU, other CPUs
were used for other jobs, so we can not guarantee exact measurement).
The second version was optimised to achieve more efficiency on the 
stage when most of non-rewritable words of length $k$ can not be extended
to non-rewritable words of length $k+1$. Concatenation of lists was 
replaced by changing the last element ``in place'', and $\proc{IsRewritableWord}$
was insered directly into the loop without a call to a separate function that,
in it turn, used {\tt ForAny}. Additionally, intersection of stabilisers
was computed in a loop which breaks if a trivial subgroup is constructed.
This and some other minor optimisations allowed to reduce the runtime to
about 15 hours. Finally, we traded space vs time and stored not only 
tuples, but also stabilisers of their elements in ${\rm Aut}(G)$. 
This permitted further speedup and reduced the runtime to be less than
ten hours. At this stage we performed six clean measurements
on a machine not running other user's jobs, and the average runtime was
9 hours and 41 minute.

The GAP code for the function {\tt RewritabilityLength} is given in the
Appendix. As you can see from the example of a GAP session below, 
the numbers of non-rewritable words 
exactly coincide with the data from \cite{Blyth-Robinson}:

\begin{verbatim}
gap> G := AlternatingGroup(5);;
gap> Exec("date");RewritabilityLength(G,10);time;Exec("date");
Wed Mar 31 11:04:39 BST 2010
Started enumeration of NRW of length 2
29 NRW of length 2 constructed
Started enumeration of NRW of length 3
1315 NRW of length 3 constructed
Started enumeration of NRW of length 4
43121 NRW of length 4 constructed
Started enumeration of NRW of length 5
528069 NRW of length 5 constructed
Started enumeration of NRW of length 6
187719 NRW of length 6 constructed
Started enumeration of NRW of length 7
1320 NRW of length 7 constructed
Started enumeration of NRW of length 8
0 NRW of length 8 constructed
8
33383583
Wed Mar 31 20:21:09 BST 2010
\end{verbatim}

Furthermore, it is easy to see that we can process independently each non-rewritable
word of the length $k$ to derive all non-rewritable words of length bigger than $k$.
Clearly, this allows parallelisation. The parallel version of the algorithm was 
implemented using the master-worker skeleton from the GAP package SCSCP \cite{scscp}.
To ensure that no data were lost, the output was modified to return the total number 
of non-rewritable words of each length summing the numbers over all parallel procedure
calls:
\begin{verbatim}
gap> RewritabilityParallel(AlternatingGroup(5),10,4);
[ 0, 29, 1315, 43121, 528069, 187719, 1320, 0 ]
\end{verbatim}
The third parameter specifies that parallel computation will be started
from the words of length four. Thus, in the beginning 1315 non-rewritable 
words of length are computed sequentially to ensure an optimal task granularity. 

The program was tested first on an 8-core Intel server, with dual quad-core Intel Xeon 
5570 2.93GHz / RAM 48 GB / CentOS Linux 5.3 with 1 master and 8 workers, and then
on a cluster consisting of three machines of the same configuration as above with 1 master 
and 24 workers. The average runtime on six measurements is given in the table below.
$$
\begin{array}{|c|r|r|c|}
\hline
\text{Number}     & \text{Runtime}    & \text{Speedup} & \text{Efficiency} \\
\text{of workers} &                   &                & \text{= speedup/nr.workers} \\
\hline
4                 & 3 \text{h } 0 \text{m } 5 \text{s} & 3.23 & 0.81       \\
8                 & 1 \text{h } 50 \text{m } 23 \text{s} & 5.26 & 0.66     \\
16                &           55 \text{m } 42 \text{s} & 10.43 & 0.65      \\
24                &           37 \text{m } 16 \text{s} & 15.59 & 0.65      \\
\hline
\end{array}
$$

To summarise, we have provided an independent verification of the result
from \cite{Blyth-Robinson}.
This may be considered as an additional motivation to find a theoretical 
proof that $A_5$ is 8-rewritable. While the original publication 20 years 
mentions the usage of PASCAL or C++ for sequential and parallel computations respectively,
now this computation has been implemented in GAP, being compatible 
with other GAP code and readable by a suitably qualified GAP user.
Parallel tools offered in the SCSCP package made it possible to 
parallelise the code within the GAP system without the necessity 
to switch to other traditional for the high-performance computing languages that
support parallelism. Note that an ongoing HPC-GAP project 
(\verb+http://www-circa.mcs.st-and.ac.uk/hpcgap.php+) is aimed to reengineer
the GAP system to provide better support for shared and distributed memory 
programming models, so in the future this example may be hopefully even better
reimplemented in a new version of the GAP system.

\newpage
\section*{Appendix: GAP Source Code for the sequential version}

\begin{verbatim}
RewritabilityLength:=function(G,limit)
local eltsG, s, A, orbsA, x, q, n, isnrw, 
      nrw, S, i, j, u, K, y, orbsK, v, tw; 
eltsG:=Filtered( G, s -> s <> () );
A:=AutomorphismGroup(G);
orbsA:=Orbits(A,G);
x:=Filtered( List(orbsA, q -> q[1] ), q -> q <> () );
x:=List( x, q -> [ [ q ], Stabilizer( A,q ) ] );
n:=1;
repeat
    n:=n+1; nrw:=[];
    Print("Started enumeration of NRW of length ", n, "\n");
    S := SymmetricGroup( n );
    S := Filtered( S, s -> s <> () );
    for i in [1..Length(x)] do   
        # Print( i, "/", Length(nrw), "\r");
        u := x[i][1]; K := x[i][2];
        if Size(K) = 1 then
            y := eltsG;        
        else
            orbsK:=Orbits(K,G);
            y:=Filtered( List(orbsK, q -> q[1] ), q -> q <> () );
        fi;    
        tw := u;
        for v in y do
            tw[n] := v;
            isnrw:=true;
            for s in S do
                if Product(tw)=Product(Permuted(tw,s)) then
                    isnrw := false; break;
                fi;    
            od;
            if isnrw then     
                Add( nrw, [ ShallowCopy(tw), 
                            Intersection( K, Stabilizer( A,v ) ) ] );
            fi;
        od;
    od;
    Print( Length(nrw), " NRW of length ", n, " constructed\n");
    if nrw=[] then return n; fi;
    x := ShallowCopy( nrw );
until n=limit;
return fail;
end;
\end{verbatim}


\begin{thebibliography}{10}

\bibitem{Blyth-Robinson}
R.~Blyth, D.~Robinson, Solution of the solubility problem for rewritable groups.  
J. London Math. Soc. (2)  41  (1990), no. 3, 438--444.

\bibitem{GAP}
The GAP~Group, \emph{GAP -- Groups, Algorithms, and Programming, 
Version 4.4.12}; 2008, \url{http://www.gap-system.org}.

\bibitem{scscp}
A.~Konovalov, S.~Linton. SCSCP --- Symbolic Computation Software Composability Protocol, 
Version 1.2; 2010, \url{http://www.cs.st-andrews.ac.uk/~alexk/scscp.htm}.

\end{thebibliography}
\end{document}